\documentstyle[a4]{article}

\input epsf
\input amssymb.sty

\def\al{\alpha}

\def\lb{\lambda}
\def\si{\sigma}

\def\ds{\oplus}
\def\map{\rightarrow}
\def\bq{\begin{equation}}
\def\eq{\end{equation}}

\def\rc{{\Bbb R}}

\def\ti{\tilde}

\def\Si{\Sigma}
\def\Ga{\Gamma}
\def\ss{\subset}

\def\cn{{\cal C}^\infty}

\begin{document}

\begin{center}
{\LARGE Contact reductions, mechanics and duality}
\vskip 0.7cm
{\large Pavol \v Severa\\}
\vskip 0.5cm
{\small \it I.H.\'E.S., Le Bois-Marie, 35, Route de Chartres \\ F-91440 
Bures-sur-Yvette, France}
\vskip 0.7cm
\end{center}
\begin{abstract}
Contact reduction is very closely related to symplectic reduction, but it allows 
 symmetries that are not manifest in Hamiltonian mechanics and moreover, 
solution of the reduced problems yields solution of the original problem without 
further integration.
\end{abstract}

\section{An informal introduction: What is the visual meaning of mechanics?}
This is a very broad question. We shall be rather modest and content ourselves 
with the answer found by Hamilton and Jacobi:
$$\epsfxsize 5cm \epsfbox{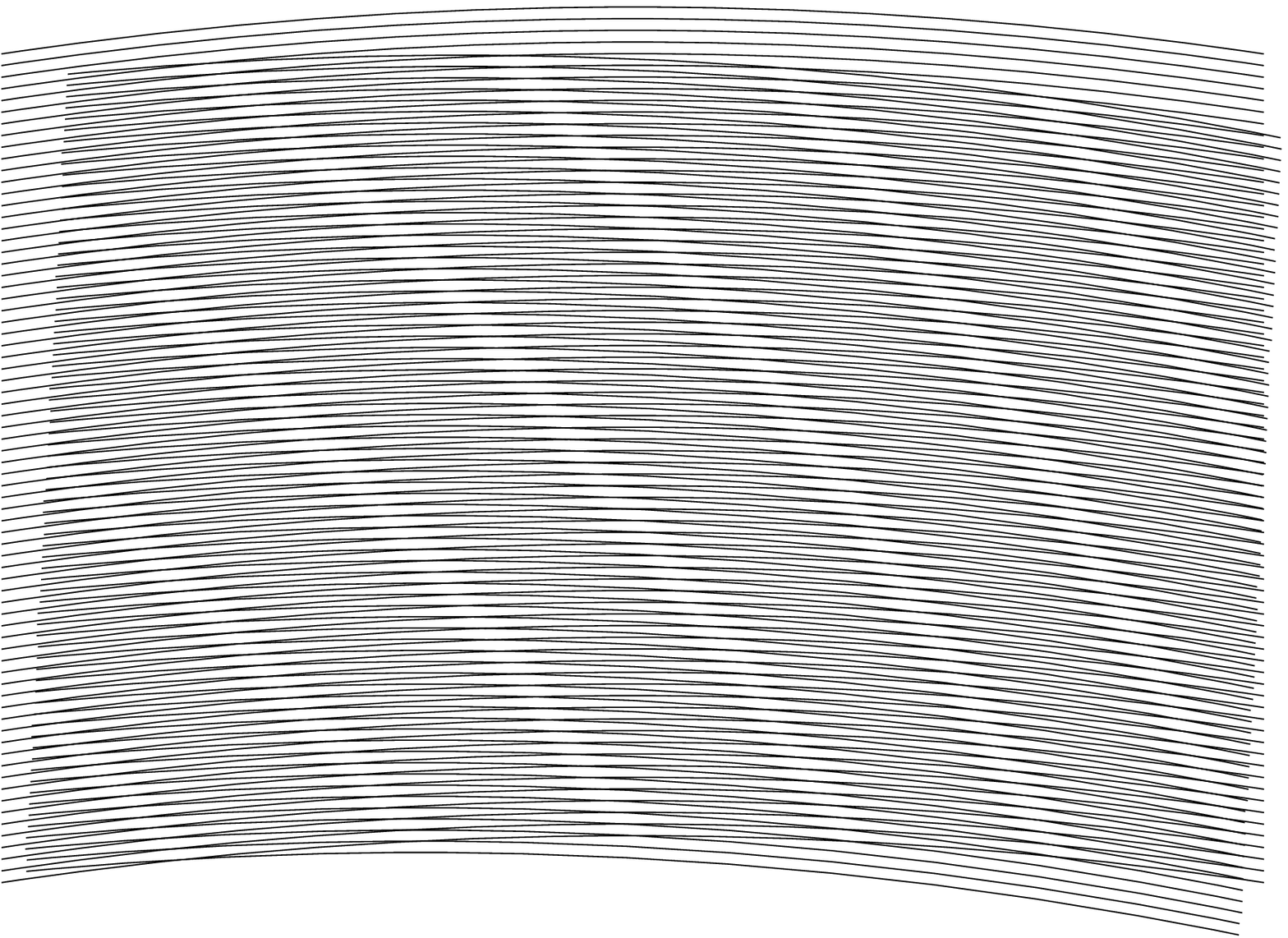}$$
According to them, mechanics studies emergence of the thick curves via 
interference of the thin curves. In higher dimensions, the thin curves are 
replaced by hypersurfaces (thick curves remain curves). The picture may 
represent the spacetime with 
worldlines of particles appearing as wave packets. A bit more abstractly, it may 
represent an extended configuration space. It may also be the ordinary space 
with rays of light. And after all, you may see similar phenomena on folded 
curtains.

The thin hypersurfaces are supposed to obey certain law (the Hamilton--Jacobi 
(or eikonal) equation). It may be stated (a bit informally) as follows:
$$\epsfxsize 5cm \epsfbox{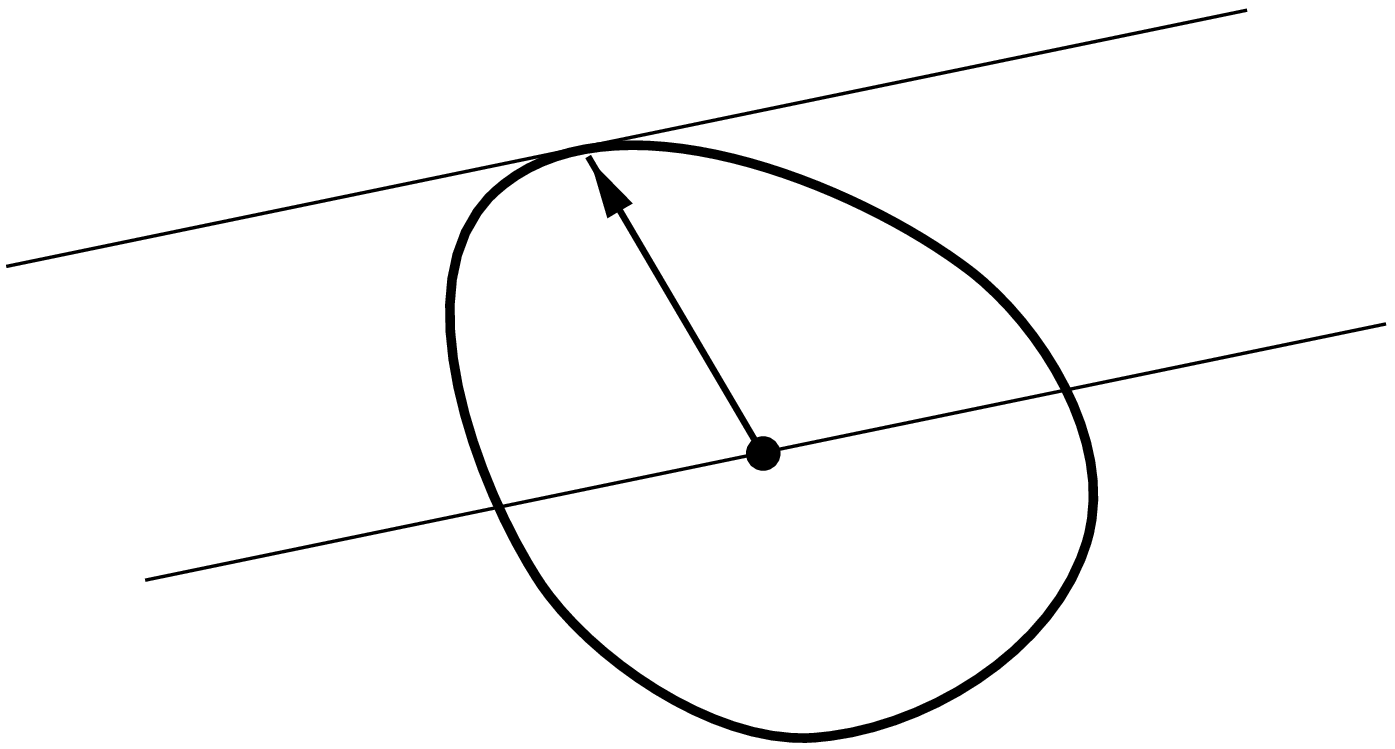}$$
Around each point we have a hypersurface (so-called wave diagram, or indicatrix) 
and if a thin hypersurface passes through the point, the next one should touch 
the diagram (the Huyghens principle). The vector represents the direction of the 
thick curve through the point. Of course, everything is understood 
infinitesimally: the picture actually represents the tangent space at the point 
and the hyperplanes are the zero- and the one-level of the differential of the 
wave phase. The wave diagram can be used to measure lengths of vectors: by 
definition, the diagram consists of vectors of the length one. Obviously, the 
thick curves are extremal with respect to this metric (this is due to the fact 
that the drawn vector points to the point of tangency). Noether theorem is 
visually obvious as well: on the first picture you can imagine that we started 
with one set of thin curves and the other one was obtained by applying an 
infinitesimal symmetry.

There are several problems that we passed in silence. The most serious one is 
non-naturalness: in Lagrangian mechanics we may add the differential of any 
function to the Lagrangian without changing anything essential. It is quite 
clear what is the corresponding fact in the wave picture: the waves are not just 
complex functions, but rather sections of a line bundle. Similarly, the phase of 
a wave is a section of the corresponding principal $U(1)$ bundle. Even if the 
bundle is trivial, only its bundle structure is natural. This modification may 
seem minute, but pictures change rather dramatically, as we have one more 
dimension: 

$$\epsfxsize 5cm \epsfbox{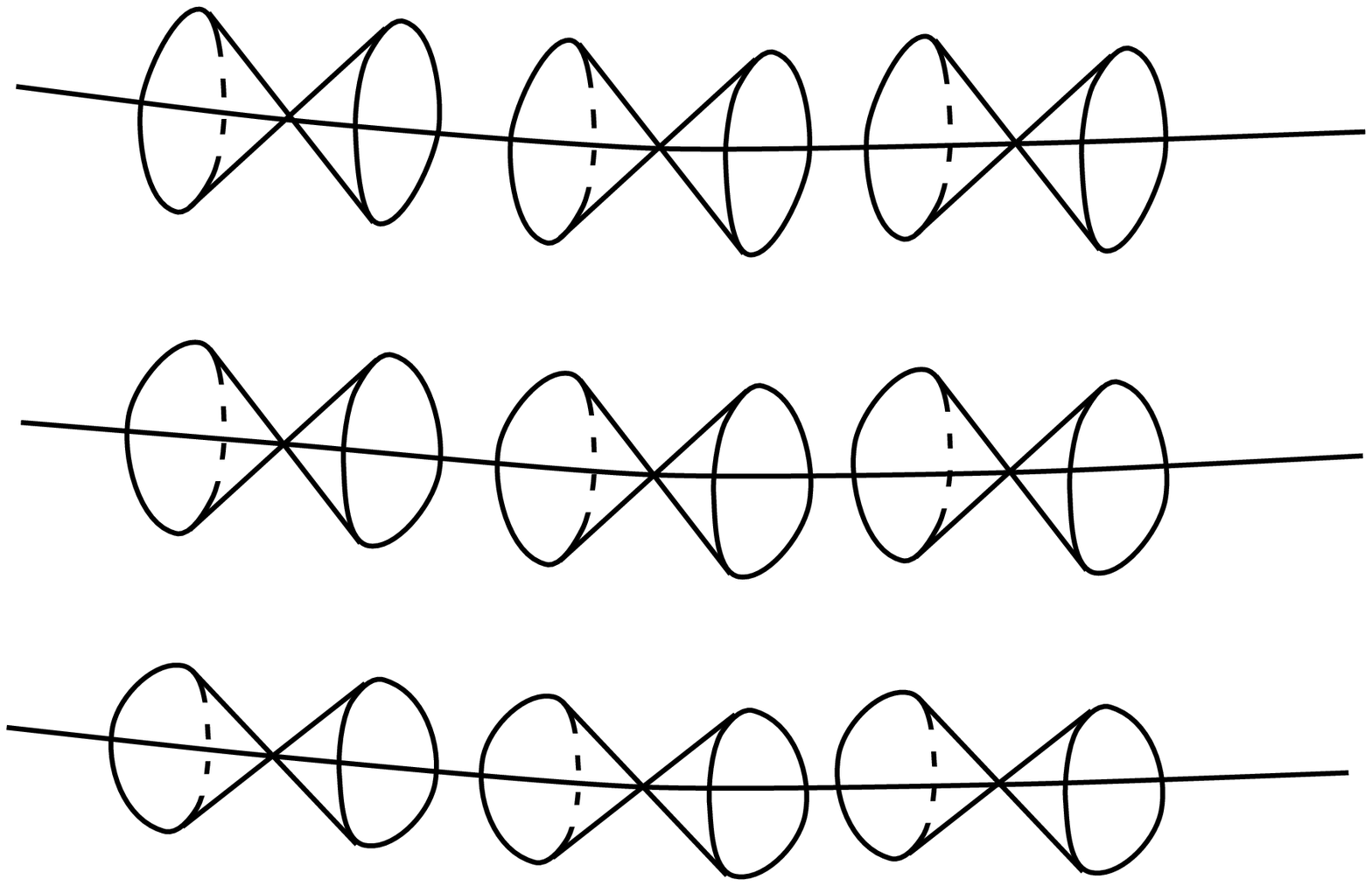}$$

In the principal $U(1)$ bundle $P$ we have an invariant field of cones. To 
obtain a Lagrangian we choose a local trivialization. Let $\al$ be the 
corresponding flat connection 1-form. We intersect the cones with the levels 
$\al=1$ and project the intersections to the base; in this way we obtain the 
wave diagrams.

Hamilton--Jacobi equation says that wave phase (a section of $P$) should be 
tangent to the cones. In other words, they are the Monge cones for the H-J 
equation. The ``worldlines'' are the (bi)characteristics, i.e. the curves along 
which the phase touches the cones.

We shall consider this picture as fundamental. We want to find characteristics 
of a field of cones and it is only a one-parameter group of symmetry that makes 
it into a variational problem. The space of characteristics is a contact 
manifold; the $U(1)$ symmetry makes it into a $U(1)$ bundle with a connection 
(ignoring some obvious problems). From this point of view, this $U(1)$ bundle 
over the phase space is more fundamental than the phase space itself. Naturality 
of this approach is advocated further in \cite{Se} and it is used in a 
substantial way in the book \cite{So}.

Now imagine that the field of cones is invariant with respect to some (local) 
Lie group $G$ containing our $U(1)$. This symmetry would be hidden from us if we 
saw only the base and not the bundle $P$ itself (unless the $U(1)$ were in its 
centre), even if we passed to Hamiltonian mechanics. We may use $G$ to reduce 
the problem (via contact reduction). The 
classical example is for two-dimensional abelian $G$. We take a vector field $v$ 
from $\frak g$, make quotient of $P$ by $v$ and find a field of cones in the 
quotient -- they are the contours of the original cones when seen in the 
direction of $v$. To find the characteristics in $P$ it is sufficient to find 
them in the reductions and to compute some derivatives. The proper setting for a 
general $G$ is contact geometry: we want to find the characteristics of a field 
of hyperplanes and we use a symmetry group $G$ to reduce the problem. Having 
solved the reductions we solve the original problem simply by computing 
derivatives. In this respect it is slightly better than symplectic reduction.

And there is one more point: we can take a different $U(1)$ in the group $G$ and 
make $P$ a bundle in this different way. The two mechanical (or variational) 
systems are thus made equivalent. This is the duality mentioned in the title.
The two systems have different phase spaces, but they share the $U(1)$-bundle 
over their phase spaces.
 
Everything is very simple and we could end here. The remaining sections should 
only provide examples and more precise definitions and make the paper completely 
selfcontained. But we should also return to our motivation and explain what is 
the high-frequency approximation in this picture. Suppose $P$ is a principal 
$U(1)$ or $\rc$ bundle and $D$ is an invariant differential operator 
$\cn(P)\map\cn(P)$. If we confine ourselves to equivariant functions with some 
weight $i/\hbar$, $D$ becomes an operator $D_\hbar$ on this associated line 
bundle. For example, if
$$D={1\over 2m}\Delta+V(x,t){\partial^2\over\partial s^2}
+{\partial^2\over\partial s\partial t}$$
and the group acts by shifting $s$, $D_\hbar$ is the  Schroedinger operator. If 
$f\in\cn(P)$ satisfies $Df=0$, each of its Fourier components satisfies $D_\hbar 
f_\hbar=0$. So this is the relation classical--quantum in our picture: all 
admissible $\hbar$'s are collected in a single equation $Df=0$ and the classical 
theory describes singularities of solutions of this wave equation.

\section{Contact geometry}

In this section we review some notions of contact geometry. Fortunately, we need 
only elementary things. Nice exposition can be found in the EMS article 
\cite{AG}.

A contact structure on a manifold $M$ is a maximally non-integrable field of 
hyperplanes $C$. For a while, let $C$ denote any subbundle of $TM$ (with 
arbitrary codimension). If $v$ is a vector field, let ${\cal L}_vC:C\map TM/C$ 
denote the infinitesimal deformation of $C$ by the flow of $v$. It may be 
defined by ${\cal L}_vC(u)=[v,u]\mbox{ mod }C$, where $u$ is extended to a 
section of $C$ in an arbitrary way. If $v$ is a section of $C$ as well, 
$[v,u]\mbox{ mod }C$ is completely local in both arguments, so that we have a 
map $\tau_C:C\wedge C\map TM/C$. If $C$ is a field of hyperplanes (a 
codimension-one subbundle), it is called maximally nonintegrable if the bilinear 
form $\tau_C$ (with values in the line-bundle $TM/C$) is everywhere regular. In 
other words, $(M,C)$ is contact iff for any infinitesimal deformation $C\map 
TM/C$ there is (unique) section of $C$ with flow generating this deformation.

\subsubsection*{Stability and characteristics}
A simple consequence of the last statement is Gray's stability theorem. If 
$(M,C(t))$ is a compact manifold with time-dependent contact structure and 
$C(t)$ is constant on some subset $X\ss M$, there is a flow of $M$ fixing the 
points of $X$ that generates $C(t)$ from $C(0)$. Compactness is used to ensure 
completness of vector fields. It is not needed if we want a local statement. For 
example, if $C(t)$ is fixed at a point (and $t$ restricted to a finite 
interval), we still have a local diffeomorphism close to the point. As a 
consequence, all contact structures in a given (necessarily odd) dimension are 
locally equivalent (Darboux theorem). 

And here is another application. A vector field $v$ on a contact manifold 
$(M,C)$ is {\it contact}, if its flow preserves $C$. For any section $w$ of 
$TM/C$ there is unique contact $v$ equal to $w$ mod $C$. To see it, take any 
$v'$ equal to $w$ mod $C$; the deformation  ${\cal L}_{v'}C$ can be removed by 
unique section of $C$. $w$ is called the {\it contact Hamiltonian} of $v$.

Contact and symplectic geometries are closely related. We shall need this 
example. Let $M\map N$ be a pricipal $U(1)$ (or $\rc$) bundle and let $M$ be 
equipped with an invariant contact structure. Suppose furthemore that the 
hyperplanes are transversal to the fibres, so that they may be interpreted as a 
connection. Its curvature is a symplectic form on $M$. A contact vector field 
commuting with the $U(1)$ action (i.e. with invariant contact Hamiltonian) is 
projected to a Hamiltonian vector field on the base:  using the connection 
1-form to trivialize $TM/C$, an invariant contact Hamiltonian becomes a function 
on the base -- the sought Hamiltonian. The somewhat mystic generation of vector 
fields by functions in symplectic geometry is here a bit more visual: 
Hamiltonian is simply the vertical part of the vector field. 
Interesting geometry arises even if $C$ is not everywhere transversal to the 
fibres. The projection of these dangerous points is easily seen to be a smooth 
hypersurface. The symplectic form diverges at the hypersurface and the 
hypersurface itself carries a contact structure (this is our first encounter 
with contact reduction).

Let us now consider a field of hyperplanes $C$ without assuming complete 
nonintegrability. A section of $C$ lying in the kernel of $\tau_C$ is called a 
{\it characteristic vector field}. The flow of a characteristic field preserves 
$C$. Suppose now that the rank of $\tau_C$ is constant. We see that the manifold 
is locally a product of a contact manifold and $\rc^k$ (the characteristic 
directions), where $k$ is the dimension of the kernel of $\tau_C$. Globally we 
have a foliation with $k$-dimensional leaves, called {\it characteristics}. If 
the foliation is actually a fibration, its base is contact.

This permits us to see the proof of Gray's theorem and the generation of contact 
fields by contact Hamiltonians from a ``space-time'' point of view. Let 
$(M,C(t))$ be as before. On $M\times\rc$ we define a hyperplane field -- we 
simply add $\partial/\partial t$ to $C(t)$. The characteristic curves are the 
worldlines of the flow. Similarly, for a given section of $TM/C$ on a contact 
$(M,C)$, we construct a hyperplane field on $M\times\rc$ in the obvious way; 
characterictics are again the worldlines of the flow.

\subsubsection*{Contact elements and 1st order PDE's}
The classical example of a contact manifold  is the space of contact elements
(i.e. hyperplanes in the tangent space) of a manifold $M$, which we denote as
$CM$. In other words, $CM$ is the projective bundle associated with $T^*M$. The 
field $C$ is given as follows: take an $x\in CM$; it 
corresponds to a hyperplane $H$ in $T_{\pi(x)}M$, where $\pi:CM\map M$ is the
natural projection. Then $(d_x\pi)^{-1}(H)$ is $C$ at $x$.

Contact geometry, in particular on $CM$, was invented by Lie to give a 
geometrical
meaning to first order PDE's and to Lagrange method
of characteristics. Suppose $E\ss CM$ is a hypersurface; it will represent
the equation. Any hypersurface $\Si\ss M$ can be lifted to $CM$: for any point
$x\in \Si$ take the hyperplane $T_x\Si$ to be a point of the lift $\ti\Si$.
$\ti\Si$ is a {\it Legendre submanifold} of $CM$, i.e. $T\ti\Si\ss C$ and 
$\ti\Si$ has the maximal dimension
(${\rm dim}\,CM=2\,{\rm dim}\,\ti\Si+1$). $\Si$
is said to solve the equation if $\ti\Si\ss E$. This has a nice interpretation
due to Monge: For any $x\in M$ we take the enveloping cone of the hyperplanes
$\pi^{-1}(x)\cap E$ in $T_xM$. In this way we obtain a field of cones in $M$.
Then $\Si$ solves the equation if it is tangent to the cones.

Lie's point of view is to forget about $M$ and to take as a solution any
Legendre submanifold contained in $E$. Such a solution may look singular in $M$
(singularities emerge upon the projection $\pi:CM\map M$; actually, many things 
classically called functions or hypersurfaces, are in fact Legendre submanifolds 
of something -- we shall meet the example of ``generating functions'' of contact 
transformations). This definition uses
only the contact structure on $CM$ and thus allows using the entire
(pseudo)group
of contact transformations.

The hyperplane field
$C$ cuts a hyperplane field $C_E$ on $E$ (there may be points where the
contact hyperplane touches $E$; generally they are isolated and we will ignore
them). The form $\tau_C$
becomes degenerate when restricted from $C$ to $C_E$;
in $E$, we have one characteristic direction everywhere.
For example, if the Monge cones coming from  $E$
are the null cones of some pseudo-Riemannian metrics on $M$ then
the projections
of the characteristics are the light-like geodesics in $M$. Legendre 
submanifolds contained in $E$ are woven from characteristics. To find them (i.e. 
to solve the equation), we have to take local quotients of $E$ by 
characteristics and pull their Legendre submanifolds back to $E$.

Hypersurfaces $E\ss CM$ often come from an equation of the type $Df=0$,
where $D:\cn (M)\map\cn (M)$ is a linear (pseudo)differential 
operator. We shall be very brief and prove nothing. Take the symbol $\si_D$ of 
$D$ (a function on $T^*M$ defined by
$D\exp(i\lb g)=(i\lb)^n\exp(i\lb g) \si_D(dg)+O(\lb^{n-1})$, $\lb\map\infty$, 
where $n$ is the
degree of $D$ and $g\in\cn(M)$).
The equation $\si_D=0$ specifies a hypersurface $E\ss CM$.
Singularities of solutions of $Df=0$ are located on ``hypersurfaces'' 
solving
the equation corresponding to $E$.

\subsubsection*{Contact reductions and homogeneous spaces}
Before proceeding to group actions and reductions, we have to describe the {\it 
contact product}. If $M_1$ and $M_2$ are contact, there is no contact structure 
on $M_1\times M_2$. The contact product $M_1\times_c M_2$ is actually a circle 
bundle over $M_1\times M_2$. In $M_1\times M_2$ there is an obvious field of 
codimension-two subspaces. We take all the hyperplanes containing these 
subspaces -- they form the manifold $M_1\times_c M_2$. Contact structure on 
$M_1\times_c M_2$ is defined similarly to the one on $CM$.

If $\phi:M_1\map M_2$ is a local diffeomorphism preserving the contact structure 
then the graph of $\phi$ in $M_1\times M_2$ can be uniquely lifted to a Legendre 
submanifold of $M_1\times_c M_2$. Vice versa, any Legendre submanifold $L\ss 
M_1\times_c M_2$ whose projection is the graph of something gives rise to such a 
$\phi$. $L$ is called (a bit improperly) the {\it generating function} of 
$\phi$. Notice that $CM\times_cCN=C(M\times N)$; the generating function of a 
$\phi:CM\map CN$ is therefore a ``hypersurface'' in $M\times N$.

Like in symplectic geometry, it is natural to call Legendre submanifolds of 
$M_1\times_cM_2\times_c\dots\times_c M_k$ contact relations; they can be 
composed under some qualifications on intersections. And similarly to $C(M\times 
M)$, Legendre submanifolds of $C\Ga$, where $\Ga$ is a Lie groupoid, can be 
composed. This represents composition of singularities in the groupoid algebra. 
The reader may like to define contact groupoids.

Now we classify contact homogeneous spaces. An invariant contact structure on 
$G/H$ can be pulled back to a left-invariand hyperplane field on $G$; its 
characteristics are cosets of $H$. Hence we are done: local $G$-homogeneous 
contact spaces are classified by hyperplanes in $\frak g$. If $\frak l \ss\frak 
g$ is such a subspace, let $\frak g_{\frak l}$ be the stabilizer ${\frak 
S}({\frak l})=\{x\in{\frak g};\;[x,{\frak l}]\ss\frak l\}$ of $\frak l$ 
intersected with $\frak l$; set $H=\exp(\frak g_{\frak l})$. For global Lie 
groups there is no reason why this $H$ should be closed. $G$-homogeneous contact 
spaces are classified by pairs $({\frak l},H)$, where $H\ss G$ is a closed 
subgroup with the Lie algebra $\frak g_{\frak l}$.

These spaces are closely related to coadjoint orbits. Choose an $\al\in\frak 
g^*$ with the kernel $\frak l$. If the orbit ${\cal O}_\al$ is homogeneous then 
$G/H={\cal O}_\al/\rc^*$ or ${\cal O}_\al/\rc_+$ (with the obvious contact 
structure coming from the symplectic form on ${\cal O}_\al$). If it is not, 
$G/H$ is a principal bundle with one-parameter group over ${\cal O}_\al$. These 
two possibilities appear according to whether ${\frak S}({\frak l})$ is 
contained in $\frak l$ or not; in the second case, $\frak S(l)/g_l$ is the 
algebra of the structure group. The constrains of globality are well visible 
here: in the first case there are none while in the second case the periods of 
the symplectic form on ${\cal O}_\al$ have to be co-mesurable (and it is also 
sufficient if $G$ is simply-conected).

Now we shall discuss contact reductions. Suppose $G$ acts on a contact space $M$ 
respecting the contact structure. Let $M_{\frak g}$ consist of the points of $M$ 
where the image of $\frak g$ is contained in $C$. If the action of $G$ is 
locally free, $M_{\frak g}$ is a submanifold, $C$ cuts a hyperplane field there 
and its characteristics are the orbits of $G$. If $G$ acts properly there (so 
that  $M_{\frak g}/G$ is a manifold),  $M//G=M_{\frak g}/G$ is the {\it contact 
reduction}; $\dim M//G=\dim M-2\dim G$. Also notice that $G$-invariant Legendre 
sumbanifolds of $M$ are contained in $M_{\frak g}$ and they are the preimages of 
Legendre submanifolds of $M//G$.

 These claims are easy to see. The contact Hamiltonians of $\frak g$ vanish at 
$M_{\frak g}$; because the action is locally free, their differentials are 
linearly independent there. It follows that $M_{\frak g}$ is a submanifold of 
codimension equal to $\dim G$ and it is nowhere touched by $C$. The vector 
fields generated by $\frak g$ are characteristic on $M_{\frak g}$ (they preserve 
the hyperplane field); for dimensional reasons they generate all the 
characteristic fields.

To reduce all the parts of $M$ (not just $M_{\frak g}$) we are sometimes forced 
to leave global geometry. Let $G$ be a (local) group acting locally freely on 
$M$. For any point $x\in M$ let ${\frak l}(x)$ consist of the vectors in $\frak 
g$ mapped into $C$ at $x$, and for any hyperplane $\frak l\ss g$ let $M_{\frak 
l}\ss M$ consist of $x$'s where ${\frak l}(x)=\frak l$. Then
$$M//_{\frak l}G=M_{\frak l}/\exp{\frak g_{\frak l}}=\left(M\times_c 
\left(G/\exp{\frak g_{\frak l}}\right)\right)//G$$
is the contact reduction at $\frak l$; if $\exp{\frak g_{\frak l}}$ is not 
closed, it makes sense only locally. For example, $CG//_{\frak l}G$ is the 
locally homogeneous space corresponding to $\frak l$. If $\frak S(l)$ is not 
contained in $\frak l$, $\frak S(l)/g_l$ is a residual one-parameter symmetry  
making $M//_{\frak l}G$ locally into a principal bundle over a symplectic space.

Let us describe two examples connected with symplectic geometry. Let $M_i\map 
N_i$ be $U(1)$ (or $\rc$) bundles with contact $M$'s and symplectic $N$'s, as 
above. Let $U(1)$ act diagonally on $M_1\times_c M_2$; then $ (M_1\times_c 
M_2)//U(1)=(M_1\times M_2)/U(1)$. The generating function (recall that it is in 
fact a Legendre submanifold) of a $U(1)$ equivariant contact map $\phi:M_1\map 
M_2$ is contained in $ (M_1\times_c M_2)_{{\frak u}(1)}$ (as it is $U(1)$ 
invariant); therefore, it is the preimage of a Legendre submanifold  $L\ss 
(M_1\times_c M_2)//U(1)$. The map $\phi$ gives rise to a symplectic map $\psi$ 
between the bases; $L$ is called a generating function of $\psi$.

Similarly, if $U(1)\times G$ acts on $M$ (the $U(1)$ is the structure group of 
the bundle), the contact reductions of $M$ give rise to the sympectic reductions 
of $N$ (a hyperplane in $\rc\times\frak g$ is the same as an element of $\frak 
g^*$).

\subsubsection*{Solving characteristic problem via contact reduction}
Finally, let us describe how contact reduction can help us in constructing the 
characteristic foliation of a hyperplane field. It is very simple: the reduction 
described above can be applied if $C$ is any hyperplane field (not necessarily 
contact). We find characteristics in the reduced manifold (this is the reduced 
problem) and we know that original characteristics are contained in their 
preimages. The problem is solved at this point: these preimages 
yield characteristics by a very simple procedure, generalizing the method of 
complete integrals. This is slightly better than symplectic reduction, since 
there solving reductions does not solve the original problem (actually, using 
the connection between symplectic and contact reduction we see that it is only 
necessary to compute an idefinite integral).

Here are the details.  Suppose $M$ is a manifold with a 
hyperplane field $C$ and that a (local) group $G$ acts on $M$ preserving $C$. 
${\frak l}(x)$ is defined as above; it is constant along the characteristics, 
because the action of $G$ on $M$ gives an action on the local quotients by 
characteristics (a contact version of Noether theorem). In other words, the 
subspaces $M_{\frak l}$ contain characteristics. We may take $M_{\frak l}/\exp 
\frak g_l$ and find characteristics there; the characteristics of $(M,C)$ are in 
their preimages. We suppose that these preimages form a foliation.

We are in this situation: we have a foliation of $M$ with 
leaves containing the characteristics and moreover the tangent spaces of the 
leaves are contained in $C$. The following is a very minute generalization of 
the method of complete integrals. Take an open subset $U\ss M$ where the 
foliation becomes a fibration $\pi:U\map B$. The hyperplanes $C$ in $U$ are 
projected to hyperplanes in $B$. This projected hyperplane is constant along 
characteristics. And vice versa -- this property specifies characteristics. The 
problem is solved.

\section{Examples}

We only consider very simple examples: right-invariant hypersurfaces $E\ss CG$ 
for various Lie groups $G$. Such a hypersurface is specified by a cone in $\frak 
g^*$ or by its dual (possibly singular) cone in $\frak g$; the field of the 
Monge cones is generated from the latter by right translation. If we choose a 
(closed, if we want to be global) one-parameter subgroup $R\ss G$, we receive a 
variational problem on $G/R$.

To find the characteristics in $E$ we use the contact reduction: we either set 
$M=E$ in the method described or we make reduction of $CG$ and look what happens 
with $E$ (these two methods are of course almost identical).

Our first example is the Euclidean group $G=SE(2)$. In this case, there are only 
two non-conjugate $\frak l$'s: either the abelian ideal of $\frak g$ or anything 
else. The first case is trivial (characteristics with this $\frak l$ are the 
straight lines in the abelian normal subgroup (or its cosets) touching the 
Monge cones). Let us reduce the second case. We find this equipment:

$$\epsfxsize 5cm \epsfbox{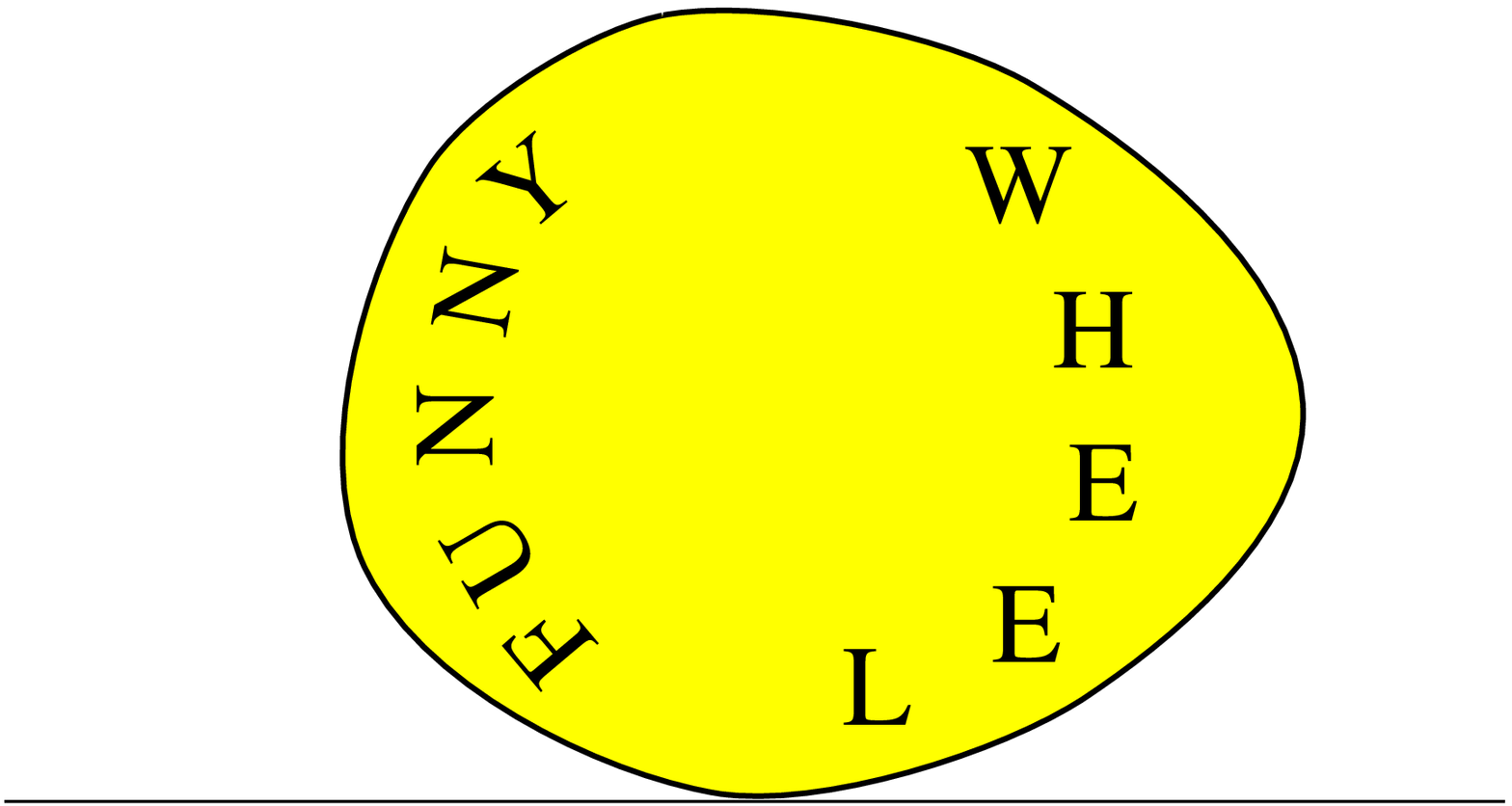}$$
The wheel is funny not just because of its shape. Instead of rolling on the road 
it remains still and it is the road that rolls around the wheel. Before being 
more explicit, let us describe the corresponding variational system in the 
Euclidean plane, i.e. in $G/R$, where $R=SO(2)$. We draw tangent vectors 
directly in the plane. To obtain the wave diagram of a point, rotate the funny 
wheel for $90^{\rm o}$ around the point. The extremals of this system are the 
trajectories swept by the points of the plane of the road; among them, the 
evolvents of the wheel (these are the trajectories of the points lying directly 
on the road). To find the trajectories it is enough to compute an indefinite 
integral (the length of arcs of the wheel). The reason behind is the residual 
one-parameter symmetry mentioned above (we shall make it explicit in a moment).

Now some details. $SE(2)$ is a semidirect product of $SO(2)$ and a vector plane 
$V_2$. The algebra $\frak g$ is (as a vector space) $\rc\ds V_2$. To find the 
Monge cone in $\frak g$, choose an origin in the plane of the funny wheel (this 
choice does not matter -- we may choose another origin of $SE(2)$), turn it for 
$90^{\rm o}$ and place it into the plane $P=\{1\}\times V_2\ss\frak g$; this 
will be 
the intersection of $P$ with the Monge cone. The adjoint action of $G$ on $P$ is 
as the defining affine action, but rotated for $90^{\rm o}$. Now we choose an 
$\frak 
l$; for definiteness, let it contain the algebra of $SO(2)$. We may visualize 
$G=SE(2)$ with the left-invariant contact structure generated by $\frak l$ 
directly in the plane $P$; in this way it becomes a double-cover of $CP$. We 
represent points of $G$ by flags (point, oriented line though the point). $1\in 
G$ is represented by $({\frak so}(2)\cap P,{\frak l}\cap P)$ (with some 
orientation of the line); acting on it by a $g\in G$ (adjoint action) we receive 
the flag of $g$. For each flag we take the orthogonal contact element -- this 
defines the contact structure. The left action of $G$ is represented by the 
adjoint action on $P$ and there is a one-parameter right action preserving the 
contact structure (a line in $V_2$) -- translating the point of the flag along 
its line.

The contact plane at $g$ touches the Monge cone if the line of its flag touches 
the wheel. Everything is very simple now. We know the projections of the 
characteristics to the quotient $G/S$, where $S$ is the one-parametric right 
symmetry, i.e. to the space of oriented lines in $P$ (this space  inherits a 
symplectic form from the contact structure on $G$). The image consists of the 
lines tangent to the wheel. We have to lift it back to $G$, moving the point 
orthogonally to the line. In other words: take an evolvent of the wheel, and 
make its points into flags by taking the lines normal to the evolvent. This 
curve in $G$ is a characteristic.

Let us mention that we may choose other $R$ to obtain a variational system -- a 
line in the vector space $V_2$. In this case, $G/R$ is the space of oriented 
lines in the Euclidean plane. We have a duality between a variational system on 
the Euclidean plane and another system on the space of its lines.

The next example is completely trivial. We choose $G$ to be the group of 
homoteties of an affine space (of arbitrary dimension) and $R$ to be the 
subgroup fixing a point. $G/R$ is the affine space. The wave diagram is one for 
all the points (warning: this is {\it not} translation-invariance). The 
extremals are simply straight lines. This example is trivial because all the 
reductions of $CG$ are one-dimensional; there is nothing to compute.

The final example is $G=SO(3)$, $R=SO(2)$. It is very similar to $SE(2)$-case. 
There is a little problem: to find the Lagrangian (or the wave diagrams) from 
the Monge cones we need a local trivialization of $G\map G/R$. In the previous 
examples we used natural global trivializations (the groups were semidirect 
products). Rather than making arbitrary choices, we use the natural (non-flat) 
connection (orthogonal to the fibres). We use it to find wave diagrams as we 
used flat connections. However, its curvature (the area form on $S^2$) has to be 
understood as a magnetic field and its potential has to be added to the 
Lagrangian (so the unnaturality is swept under this rug).

The system looks as follows. Choose a Monge cone in $\frak g$. Identify $G/R$ 
with the unit sphere in $\frak g$. To find the wave diagram at a point of the 
sphere, take the tangent plane, intersect with the cone and rotate for $90^{\rm 
o}$. 
This time all the $\frak l$'s are mutually conjugated. The funny wheel idea 
works as before (the wheel is the intersection of the sphere with the Monge 
cone), just instead of planes and lines we have spheres and great circles. The 
problem is reduced to computing lengths of arcs of the wheel. Characteristics 
are closed if the perimeter of the wheel is a rational multiple of $\pi$. 

There is an alternative description of the way the road sphere moves on the 
wheel sphere. Intersect the dual cone with the sphere and consider the paralel 
transport along the curve; it extends naturally to a motion of the sphere. This 
is familiar from the motion of rigid body \cite{LL}. It is no accident: the body 
is described by a field of cones in $SO(3)\times\rc^2$ (one $\rc$ is the time 
and the other is the action); the reductions of this group are as those of 
$SO(3)$.

\vskip 5mm
e-mail: severa@ihes.fr


\begin{thebibliography}{9}
\bibitem{AG}V.I. Arnold, A.B. Givental, {\it Symplectic geometry}, in 
Encyclopaedia of Mathemathical Sciences 4, Springer 1990
\bibitem{LL}L.D. Landau, E.M. Lifshitz, {\it Mechanics}, Nauka, Moscow 1973
\bibitem{Se}P. \v Severa, {\it Contact geometry in Lagrangean mechanics}, 
J.Geom.Phys.(29)3 (1999) pp. 235-242
\bibitem{So}J.-M.Souriau, {\it Structure of dynamical systems}, Birkh\"auser, 
Boston 1997
\end{thebibliography}
\end{document}